\documentclass[11pt]{article}
\usepackage{amsfonts}
\usepackage{amsmath}
\usepackage{amssymb}
\begin{document}

\baselineskip 18pt
\title{\Large\bf On the generalization of Golomb's conjecture}
\author{Chaohua Jia}
\date{}
\maketitle {\small \noindent {\bf Abstract.} Let $p$ be a
sufficiently large prime number, $r$ be any given positive integer.
Suppose that $a_1,\,\dots,\,a_r$ are pairwise distinct and not zero
modulo $p$. Let $N(a_1,\,\dots,\,a_r;\,p)$ denote the number of
$\alpha_1,\,\dots,\,\alpha_r,\,\beta$, which are primitive roots
modulo $p$, such that $\alpha_1+\beta\equiv
a_1,\,\dots,\,\alpha_r+\beta\equiv a_r\,({\rm mod}\,p).$ In the
first version of this paper, we proved an asymptotic formula for
$N(a_1,\,\dots,\,a_r;\,p)$ so that we could answer an open problem
of Wenpeng Zhang and Tingting Wang. But we found that our result had
been included in a paper of L. Carlitz in 1956, which is explained
in the additional remark below.}

\vskip.3in
\noindent{\bf 1. Introduction}

Let $p$ be an odd prime number, {\bf A}(p) denote the set of all
primitive roots modulo $p$. The Golomb's conjecture[3] states that
there exist $\alpha,\,\beta\in{\bf A}(p)$ such that
$$
\alpha+\beta\equiv 1\,({\rm mod}\,p).
$$

There is extensive study on the Golomb's conjecture. Juping Wang[5]
basically solved the Golomb's conjecture for finite fields. Some
other works can be found in [6], [2].

Recently, Wenpeng Zhang and Tingting Wang[7] proved the following
theorem to generalize the Golomb's conjecture.

{\bf Theorem}(Wenpeng Zhang, Tingting Wang). \emph{For the positive
integers $a,\,b(1\leq a\ne b\leq p-1)$, let $N(a,\,b;\,p)$ denote
the number of $\alpha,\,\beta,\,\gamma\in{\bf A}(p)$ such that}
$$
\alpha+\gamma\equiv a,\qquad\quad\beta+\gamma\equiv b\,({\rm
mod}\,p).
$$
\emph{Then one has an asymptotic formula}
$$
N(a,\,b;\,p)={\varphi^3(p-1)\over p^2}+O(p^{{1\over 2}+
\varepsilon}),
$$
where $\varphi(n)$ is the Euler totient function, $\varepsilon$ is
any small positive constant.

They[7] also put forward the following open problem.

{\bf Problem}(Wenpeng Zhang, Tingting Wang). \emph{Let ${\Bbb F}_p$
be the field with $p$ elements. For pairwise distinct non-zero
elements $a,\,b,\,c\in{\Bbb F}_p$, if there exist four elements
$\alpha,\,\beta,\, \gamma,\,\delta\in{\bf A}(p)$ such that}
$$
\alpha+\delta\equiv a,\quad \beta+\delta\equiv b,\quad
\gamma+\delta\equiv c\,({\rm mod}\,p)?
$$

In this paper, we shall answer this problem for sufficiently large
prime $p$ and further generalize the Golomb's conjecture as follows.

{\bf Theorem}. \emph{Let $p$ be a sufficiently large prime number,
$r$ be any given positive integer. Suppose that
$a_1,\,\dots,\,a_r\in{\Bbb F}_p$ are pairwise distinct and not zero
modulo $p$. Let $N(a_1,\,\dots,\,a_r;\,p)$ denote the number of
$\alpha_1,\,\dots,\,\alpha_r,\,\beta\in{\bf A}(p)$ such that}
$$
\alpha_1+\beta\equiv a_1,\,\quad\dots,\,\quad\alpha_r+\beta\equiv
a_r\,({\rm mod}\,p).
$$
\emph{Then we have an asymptotic formula}
$$
N(a_1,\,\dots,\,a_r;\,p)={\varphi^{r+1}(p-1)\over p^r}+O(p^{{1\over
2}+\varepsilon}),
$$
\emph{where $\varphi(n)$ is the Euler totient function,
$\varepsilon$ is any small positive constant.}

When $r=2$, Theorem is the theorem of Wenpeng Zhang and Tingting
Wang. If $r=3$, by the fact that
$$
\varphi(p-1)\gg{p\over \log\log p},
$$
then
$$
N(a_1,\,a_2,\,a_3;\,p)\gg{p\over \log\log^3 p},
$$
which answers the problem of Wenpeng Zhang and Tingting Wang for
sufficiently large prime $p$.

\vskip.3in
\noindent{\bf 2. Two lemmas}

{\bf Lemma 1}. \emph{For any integer $a((a,\,p)=1)$, we have}
\begin{align*}
&\ \,{\varphi(p-1)\over p-1}\sum_{k|\,p-1}{\mu(k)\over \varphi(k)}
\sum_{\substack{d=1\\ (d,\,k)=1}}^k e\Bigl({d\,{\rm ind}\,a\over
k}\Bigr)\\
&=
\begin{cases}
1,\qquad\quad \emph{if}\ \ a\in{\bf A}(p),\\
0,\qquad\quad \emph{if}\ \ a\not\in{\bf A}(p),
\end{cases}
\end{align*}
\emph{where $e(x)=e^{2\pi ix},\,{\rm ind}\,a$ is the index of $a$
modulo $p$.}

One could see Proposition 2.2 in [4].

{\bf Lemma 2}. \emph{Let $\chi_1,\,\dots,\,\chi_s$ be Dirichlet
characters molulo $p$, at least one of which is non-principal. Let
$b_1,\,\dots,\,b_s\in{\Bbb F}_p$ be pairwise distinct. Then we have}
$$
\Bigl|\sum_{x=1}^p\chi_1(x+b_1)\cdots\chi_s(x+b_s)\Bigr|\leq
sp^{1\over 2}.
$$

One could see Lemma 17 in [1].

\vskip.3in
\noindent{\bf 3. Proof of Theorem}

For any integers $k,\,d(1\leq k\leq p-1,\,k|\,p-1,\,1\leq d\leq
k,\,(d,\,k)=1)$, we write
$$
e\Bigl({d\,{\rm ind}\,a\over k}\Bigr)=\chi_{d,\,k}(a),
$$
and for $p|\,a$, write
$$
\chi_{d,\,k}(a)=0.
$$

It is easy to see that $\chi_{d,\,k}(a)$ is a Dirichlet character
modulo $p$ and that $\chi_{d,\,k}(a)$ is the principal character if
and only if $k=1$.

By Lemma 1, we have
\begin{align*}
&\ \,N(a_1,\,\dots,\,a_r;\,p)\\
&=\sum_{\substack{\beta=1\\ \beta\in{\bf A}(p)}}^{p-1} \sum_
{\substack{\alpha_1=1\\ \alpha_1\in{\bf A}(p)\\
\alpha_1+\beta\equiv a_1\,({\rm
mod}\,p)}}^{p-1}\cdots\sum_{\substack{\alpha_r=1\\ \alpha_r\in{\bf
A}(p)\\ \alpha_r+\beta\equiv a_r\,({\rm mod}\,p)}}^{p-1}\\
&=\sum_{\beta=1}^{p-1}\Bigl({\varphi(p-1)\over
p-1}\sum_{k|\,p-1}{\mu(k)\over \varphi(k)}\sum_{\substack{d=1\\
(d,\,k)=1}}^k\chi_{d,\,k}(\beta)\Bigr)\cdot
\end{align*}
\begin{align*}
&\cdot\sum_{\substack{\alpha_1=1\\ \alpha_1+\beta\equiv a_1\,({\rm
mod}\,p)}}^{p-1}\Bigl({\varphi(p-1)\over
p-1}\sum_{k_1|\,p-1}{\mu(k_1)\over \varphi(k_1)}\sum_{\substack{d_1=1\\
(d_1,\,k_1)=1}}^{k_1}\chi_{d_1,\,k_1}(\alpha_1)\Bigr)\cdot\\
&\cdots\sum_{\substack{\alpha_r=1\\ \alpha_r+\beta\equiv a_r\,({\rm
mod}\,p)}}^{p-1}\Bigl({\varphi(p-1)\over
p-1}\sum_{k_r|\,p-1}{\mu(k_r)\over \varphi(k_r)}\sum_{\substack{d_r=1\\
(d_r,\,k_r)=1}}^{k_r}\chi_{d_r,\,k_r}(\alpha_r)\Bigr)\\
&={\varphi^{r+1}(p-1)\over
(p-1)^{r+1}}\Bigl(\sum_{k|\,p-1}{\mu(k)\over
\varphi(k)}\sum_{\substack{d=1\\ (d,\,k)=1}}^k\sum_{\beta=1}^p
\chi_{d,\,k}(\beta)\Bigr)\cdot\\
&\cdot\Bigl(\sum_{k_1|\,p-1}{\mu(k_1)\over
\varphi(k_1)}\sum_{\substack
{d_1=1\\ (d_1,\,k_1)=1}}^{k_1}\sum_{\substack{\alpha_1=1\\
\alpha_1\equiv a_1-\beta\,({\rm
mod}\,p)}}^p\chi_{d_1,\,k_1}(\alpha_1)\Bigr)\cdot\\
&\cdots\Bigl(\sum_{k_r|\,p-1}{\mu(k_r)\over
\varphi(k_r)}\sum_{\substack
{d_r=1\\ (d_r,\,k_r)=1}}^{k_r}\sum_{\substack{\alpha_r=1\\
\alpha_r\equiv a_r-\beta\,({\rm
mod}\,p)}}^p\chi_{d_r,\,k_r}(\alpha_r)\Bigr)\\
&={\varphi^{r+1}(p-1)\over
(p-1)^{r+1}}\Bigl(\sum_{k|\,p-1}{\mu(k)\over
\varphi(k)}\sum_{\substack{d=1\\ (d,\,k)=1}}^k\sum_{\beta=1}^p
\chi_{d,\,k}(\beta)\Bigr)\cdot\\
&\cdot\Bigl(\sum_{k_1|\,p-1}{\mu(k_1)\over
\varphi(k_1)}\sum_{\substack {d_1=1\\
(d_1,\,k_1)=1}}^{k_1}\chi_{d_1,\,k_1}(a_1-\beta)\Bigr)\cdot\\
&\cdots\Bigl(\sum_{k_r|\,p-1}{\mu(k_r)\over\varphi(k_r)}\sum_
{\substack{d_r=1\\ (d_r,\,k_r)=1}}^{k_r}\chi_{d_r,\,k_r}(a_r-
\beta)\Bigr)\\
&={\varphi^{r+1}(p-1)\over (p-1)^{r+1}}\sum_{k|\,p-1}{\mu(k)\over
\varphi(k)}\sum_{\substack{d=1\\
(d,\,k)=1}}^k\sum_{k_1|\,p-1}{\mu(k_1)\over
\varphi(k_1)}\sum_{\substack
{d_1=1\\ (d_1,\,k_1)=1}}^{k_1}\\
&\cdots\sum_{k_r|\,p-1}{\mu(k_r)\over \varphi(k_r)}\sum_{\substack
{d_r=1\\ (d_r,\,k_r)=1}}^{k_r}\sum_{\beta=1}^p
\chi_{d,\,k}(\beta)\chi_{d_1,\,k_1}(a_1-\beta)\cdots\chi_{d_r,\,
k_r}(a_r-\beta)\\
&={\sum}_1+{\sum}_2,
\end{align*}
where ${\sum}_1$ is the sum in which
$(k,\,k_1,\,\dots,\,k_r)=(1,\,1,\,\dots,\,1)$, ${\sum}_2$ is the sum
in which $(k,\,k_1,\,\dots,\,k_r)\ne(1,\,1,\,\dots,\,1)$.

The application of Lemma 2 yields that in ${\sum}_2$, we have
\begin{align*}
&\ \,\Bigl|\sum_{\beta=1}^p
\chi_{d,\,k}(\beta)\chi_{d_1,\,k_1}(a_1-\beta)\cdots\chi_{d_r,\,
k_r}(a_r-\beta)\Bigr|\\
&=\Bigl|\sum_{\beta=1}^p
\chi_{d,\,k}(\beta)\chi_{d_1,\,k_1}(\beta-a_1)\cdots\chi_{d_r,\,
k_r}(\beta-a_r)\Bigr|\\
&\leq(r+1)p^{1\over 2}.
\end{align*}
Hence,
\begin{align*}
{\sum}_2&\ll{\varphi^{r+1}(p-1)\over (p-1)^{r+1}}
\sum_{k|\,p-1}{|\mu(k)|\over\varphi(k)}\sum_{\substack{d=1\\
(d,\,k)=1}}^k\sum_{k_1|\,p-1}{|\mu(k_1)|\over
\varphi(k_1)}\sum_{\substack{d_1=1\\ (d_1,\,k_1)=1}}^{k_1}\\
&\cdots\sum_{k_r|\,p-1}{|\mu(k_r)|\over \varphi(k_r)}\sum_{\substack
{d_r=1\\ (d_r,\,k_r)=1}}^{k_r}p^{1\over 2}\\
&\ll p^{1\over 2}\sum_{k|\,p-1}|\mu(k)|\sum_{k_1|\,p-1}|\mu(k_1)|
\cdots\sum_{k_r|\,p-1}|\mu(k_r)|\\
&\ll p^{1\over 2}2^{(r+1)\omega(p-1)}\\
&\ll p^{{1\over 2}+\varepsilon},
\end{align*}
where $\omega(n)$ denotes the number of different prime factors of
$n$.

We have
\begin{align*}
{\sum}_1&={\varphi^{r+1}(p-1)\over (p-1)^{r+1}}\sum_{\beta=1}^p
\chi_0(\beta)\chi_0(a_1-\beta)\cdots\chi_0(a_r-\beta)\\
&={\varphi^{r+1}(p-1)\over (p-1)^{r+1}}(p+O(1))\\
&={\varphi^{r+1}(p-1)\over p^r}+O(1),
\end{align*}
where $\chi_0$ denotes the principal character modulo $p$.
Therefore,
$$
N(a_1,\,\dots,\,a_r;\,p)={\varphi^{r+1}(p-1)\over p^r}+O(p^{{1\over
2}+\varepsilon}).
$$

So far the proof of Theorem is finished.

\vskip.3in
\noindent{\bf 4. Additional remark}

After the first version of this paper was published in arXiv, we
found that the result in Theorem had been included in Theorem 1 of
the paper of L. Carlitz (\emph{Sets of primitive roots}, Compositio
Math., \textbf{13}(1956), 65-70). L. Carlitz proved the results in
quite general situation. If taking
$$
f_1(x)=x,\quad f_2(x)=-x+a_1,\quad\dots,\quad f_{r+1}(x)=-x+a_r
$$
and $e_i=p-1$, this is the result in our Theorem. In the end of
Section 2 of his paper, L. Carlitz pointed out that by Weil's bound,
error term could be $O(p^{{1\over 2}+\varepsilon})$.

\vskip.3in
\noindent{\bf Acknowledgements}

This work is supported by the National Natural Science Foundation of
China (Grant No. 11771424). The author would like to thank Professor
Wenpeng Zhang for his fascinating talk to introduce his new joint
work with Tingting Wang, which was given in the conference ``Number
Theory on the Cloud'' on April 18th, 2020.

The author would like to thank Professor Ke Gong for telling him the
paper of L. Carlitz. He also would like to thank Dr. Tim Trudgian
for comments and information.


\vskip.6in

Chaohua Jia

Institute of Mathematics, Academy of Mathematics and Systems
Science, Chinese Academy of Sciences, Beijing 100190, P. R. China

Hua Loo-Keng Key Laboratory of Mathematics, Chinese Academy of
Sciences, Beijing 100190, P. R. China

School of Mathematical Sciences, University of Chinese Academy of
Sciences, Beijng 100049, P. R. China

E-mail: {\tt jiach@math.ac.cn}

\end{document}